\documentclass[11pt]{article}

\usepackage{graphicx,enumitem}
\usepackage{graphicx,psfrag,epsf}
{%
\centering\addtocounter{figure}{1}
\begin{enumerate}[%
itemsep=2pt,parsep=0em,
label={(\alph*)},
ref={\thefigure.(\alph*)}
]}%
{\end{enumerate}\addtocounter{figure}{-1}}

\usepackage[colorlinks=true,linkcolor=blue]{hyperref}%
\usepackage{graphicx}
\hypersetup{urlcolor=blue, citecolor=blue, colorlinks=true, linkcolor=blue} 
 
\usepackage[margin=1in]{geometry}

\usepackage[utf8]{inputenc}
\usepackage{amsmath}
\usepackage[english]{babel}
\usepackage{amssymb}
\usepackage{amsfonts}
\usepackage{bm}
\usepackage{graphicx}
\usepackage{amsmath}
\usepackage{amsthm}
\usepackage{xcolor}
\usepackage{parskip}
\usepackage{multicol}
\usepackage{float}

\usepackage{verbatim}

\floatstyle{plaintop}
\restylefloat{table}
\usepackage[tableposition=top]{caption}

\usepackage{rotating}
\usepackage{subcaption}
\usepackage{multirow}

\usepackage{chngcntr}
\usepackage{apptools}
\AtAppendix{\counterwithin{lem}{subsection}}
\AtAppendix{\counterwithin{proposition}{subsection}}
\AtAppendix{\counterwithin{equation}{subsection}}

\usepackage{authblk}

 \theoremstyle{definition}

\newtheorem{cor}{Corollary}[section]

\newtheorem{proposition}{Proposition}[section]

\newtheorem{prueba}{Proof}[section]

\newtheorem{remark}{Remark}[section]

\usepackage{booktabs}
\usepackage{animate}

\usepackage{natbib}
\hypersetup{urlcolor=blue, citecolor=blue, colorlinks=true, linkcolor=blue} 

\numberwithin{equation}{section}
\numberwithin{figure}{section}
\numberwithin{table}{section}

\renewcommand{\qed}{\hfill\blacksquare}

\usepackage{blindtext}
\usepackage{lineno}

\begin{document}

\def\spacingset#1{\renewcommand{\baselinestretch}%
{#1}\small\normalsize} \spacingset{1}

\title{Versatile Parametric Classes of Covariance Functions that Interlace Anisotropies and Hole Effects}

\author[1]{Alfredo Alegr\'ia\footnote{Corresponding author. Email: alfredo.alegria@usm.cl}}

\affil[1]{Departamento de Matem\'atica, Universidad T{\'e}cnica Federico Santa Mar{\'i}a, Chile}

\author[2,3]{Xavier Emery}
\affil[2]{Department of Mining Engineering\\ University of Chile\\ Santiago, Chile.}
\affil[3]{Advanced Mining Technology Center\\ University of Chile\\ Santiago, Chile.}

\maketitle

\begin{abstract} 
\noindent Covariance functions are a fundamental tool for modeling the dependence structure of spatial processes. This work investigates novel constructions for covariance functions that enable the integration of anisotropies and hole effects in complex and versatile ways, having the potential to provide more accurate representations of dependence structures arising with real-world data. We show that these constructions extend widely used covariance models, including the Matérn, Cauchy, compactly-supported hypergeometric and cardinal sine models. We apply our results to a geophysical data set from a rock-carbonate aquifer and demonstrate that the proposed models yield more accurate predictions at unsampled locations compared to basic covariance models. \\

\end{abstract}

\section{Introduction}
\label{sec:intro}

Data indexed by spatial (hereafter, Euclidean) coordinates arise in many disciplines of the natural sciences, including climatology \citep{sang2011covariance}, oceanography \citep{wikle2013modern}, environment \citep{rodrigues2015non}, ecology \citep{finley2011hierarchical} and geosciences \citep{davis2002}. Statistical and geostatistical models often assume the observed data to be a realization of a Gaussian random field, with the covariance function being the fundamental ingredient to capture the spatial dependence \citep{chiles2009geostatistics}, to understand the underlying spatial patterns and to make reliable predictions. 

Currently, there is a fairly extensive catalog of parametric families of stationary covariance functions that allow modeling a large number of patterns appearing in real situations, such as long-memory, hole effects, periodicities, degree of mean square differentiability, anisotropies, among others. 
Classical textbooks, such as \cite{gaetan2010spatial} and \cite{chiles2009geostatistics}, provide extensive insights into the wide range of available models. While existing models can handle many common patterns found in real data sets, some data sets may present complex combinations of features that require the development of new specialized models. In particular, anisotropies and hole effects are two common properties that can manifest on the covariance structure of data. Anisotropy refers to the directional dependence of spatial data, where the level of association varies across different directions. We refer the reader to \cite{allard2016anisotropy} and \cite{koch2020computationally} for discussions on various types of anisotropy. Hole effects, on the other hand, refer to the occurrence of negative covariance values at large distances, which can be attributed to the structured occurrence of high (low) values of a georeferenced variable surrounded by low (high) values of this variable \citep{chiles2009geostatistics}. 

Although some basic constructions that incorporate both anisotropy and hole effects can be designed easily (some examples are provided in Section \ref{sec:background}), more complex and sophisticated relationships may be required in practice. Our focus is on covariance models that feature both amenable expressions and interpretable parameters, and that are capable of achieving  negative values of varying intensities depending on the spatial orientation. In particular, some models could display negative values only along specific spatial directions. 
We are motivated to study this type of models in order to have a flexible framework capable of capturing intricate dependence patterns present in real-world data, and enable more robust inference and prediction.

To accomplish this goal, we begin by examining the conditions under which the difference between two geometrically anisotropic stationary covariance functions is valid. In a purely isotropic setting, \cite{ma2005linear}, \cite{buhmann2020}, \cite{faouzi2020} and  \cite{posa} utilized this methodology for constructing models with hole effects. Our findings thus expand upon these works by considering an anisotropic setting. Furthermore, we investigate an approach based on the difference between a merely isotropic model and the average of shifted isotropic models. The shift direction is a critical element of this formulation as it indicates the primary direction where the hole effect occurs. 
In addition, we study a construction that involves directional derivatives of a spatial process; thus, a significant hole effect is expected in a predominant direction (directional derivative's sign can amplify the transitions between high and low values). We also investigate how the aforementioned constructions can be coupled with popular existing covariance models, such as the Matérn, Cauchy, compactly-supported hypergeometric and cardinal sine, to generalize these models to more versatile parametric functions. 

The practical implications of this work will be explored through an application to a geophysical dataset. 
Our analysis will reveal that the proposed models lead to substantially improved predictions at unsampled locations in comparison with basic covariance models.

The article is organized as follows. Section \ref{sec:background} contains preliminary material on stationary spatial random fields, covariance functions and basic models that combine anisotropies and hole effects. Section \ref{sec:results} proposes general methodologies to construct models merging anisotropies and hole effects in a nontrivial manner. Section \ref{parametric} offers explicit parametric families that use Matérn, Cauchy, compactly-supported hypergeometric and cardinal sine models as a starting point. In Section \ref{sec:data}, our findings are applied to a real data set. Section \ref{sec:conclusions} presents conclusions and outlines potential avenues for future research.

\section{Preliminaries}
\label{sec:background}

Let $d$ be a positive integer and $\{Z(\bm{x}):\bm{x}\in\mathbb{R}^d\}$ be a second-order zero-mean random field. The covariance function of such a random field is the mapping $K: \mathbb{R}^d \times \mathbb{R}^d \rightarrow \mathbb{R}$ defined as $K(\bm{x},\bm{x}') = {\rm cov}\left[Z(\bm{x}),Z(\bm{x}')\right]$. 
This is a positive semidefinite function, i.e., for all $n\in\mathbb{N}$, $v_1,\hdots,v_n\in\mathbb{R}$ and $\bm{x}_1,\hdots,\bm{x}_n \in \mathbb{R}^d$, 
\begin{equation*}
    \sum_{i,j=1}^n v_i v_j K(\bm{x}_i,\bm{x}_j) \geq 0.
\end{equation*}
 The mapping $K$ is said to be stationary if there exists a function $C:\mathbb{R}^d\rightarrow \mathbb{R}$ such that
$K(\bm{x},\bm{x}') = C(\bm{x}-\bm{x}')$, for all $\bm{x},\bm{x}'\in\mathbb{R}^d$. By abuse of language, $C$ will be referred to as a stationary covariance function and we will say that $C$ is positive semidefinite. Bochner's theorem (see, e.g., page 24 of \citealp{stein-book}) provides a useful characterization of these mappings under an assumption of continuity: $C$ is a continuous stationary covariance function if and only if it can be written as
\begin{equation}
\label{bochner}
    C(\bm{h}) = \int_{\mathbb{R}^d} \exp\left( \imath \bm{h}^\top \bm{\omega} \right) F(\text{d}\bm{\omega}), \qquad \bm{h}\in\mathbb{R}^d,
\end{equation}
for some nonnegative finite measure $F$ (called spectral measure), with $\imath$ standing for the imaginary unit. If $F$ is absolutely continuous with respect to the Lebesgue measure, which happens if $C$ is absolutely integrable, then $F(\text{d}\bm{\omega}) = f(\bm{\omega}) \text{d}\bm{\omega}$, for some function $f:\mathbb{R}^d\rightarrow \mathbb{R}$ known as the spectral density. In such a case, Fourier inversion yields
\begin{equation}
    \label{fourier_inversion}
    f(\bm{\omega}) = \frac{1}{(2\pi)^d} \int_{\mathbb{R}^d} \exp\left( -\imath \bm{\omega}^\top \bm{h}\right)C(\bm{h})
\text{d}\bm{h}, \qquad \bm{\omega}\in\mathbb{R}^d.
\end{equation}
A stationary covariance function is said to be isotropic if there exists a function $\varphi:[0,\infty) \rightarrow \mathbb{R}$ such that
$K(\bm{x},\bm{x}') = \varphi(\|\bm{x}-\bm{x}'\|)$,  for all $\bm{x},\bm{x}'\in\mathbb{R}^d$. The function $\varphi$ is referred to as the isotropic part of $K$. We denote $\Phi_d$ the set of continuous functions $\varphi$ that are the isotropic part of some positive semidefinite function in $\mathbb{R}^d \times \mathbb{R}^d$. Every member of $\Phi_d$, for $d\geq 2$, can be written as the Hankel transform of order $(d-2)/2$ of a nondecreasing bounded measure $G_d$ on $[0,\infty)$ \citep{schoenberg1938metric}, i.e.,
\begin{equation}
    \label{hankel}
    \varphi(h) =  \int_0^\infty \Omega_d(hu) \text{d}G_d(u), \qquad h \geq 0,
\end{equation}
where $\Omega_d(s) = 2^{(d-2)/2} \Gamma(d/2) s^{-(d-2)/2} J_{(d-2)/2}(s)$, with $\Gamma$ standing for the gamma function and $J_\nu$ for the Bessel function of the first kind of order $\nu$ \citep{Olver}. If the spectral measure $F$ is absolutely continuous with respect to the Lebesgue measure, then so is $G_d$ and one has
\begin{equation}
    \label{fourier_inversion2}
    \varphi(h) = {(2\pi)^{d/2}} h^{(2-d)/2} \int_0^{\infty} J_{(d-2)/2}(u h) f_d(u) u^{d/2} \text{d}u, \qquad h \geq 0,
\end{equation}
and
\begin{equation}
    \label{fourier_inversion3}
    f_d(u) = \frac{1}{(2\pi)^{d/2}} u^{(2-d)/2} \int_0^{\infty} J_{(d-2)/2}(u h) \varphi(h) h^{d/2} \text{d}h, \qquad u \geq 0,
\end{equation}
where $f_d$ is the radial part of $f$ and will be referred to as the $d$-radial spectral density of $\varphi$ (note that the expression of this radial density depends on the space dimension $d$): $f(\bm{\omega}) = f_d(\| \bm{\omega} \|)$ for all $\bm{\omega} \in \mathbb{R}^d$.  

As described in the introduction, the isotropic part of an isotropic covariance function $\varphi$ can attain negative values at large distances, which is commonly referred to as a hole effect. For simplicity, suppose that $\int_0^\infty \text{d}G_d(u) = 1$, then one has the following lower bound for the members of $\Phi_d$: 
$$\varphi(h) \geq \inf_{s\geq 0} \Omega_d(s).$$
When $d=2$ and $d=3$, this lower bound is $-0.403$ and $-0.218$, respectively \citep{stein-book}. As the spatial dimension $d$ approaches infinity, the lower bound of the isotropic covariance function tends to zero, indicating that an isotropic hole effect becomes negligible with large spatial dimensions.

In the following sections, we aim to investigate parametric covariance models that interlace anisotropy and hole effect. Note that some elementary constructions can be developed: 
\begin{itemize}
\item Suppose that $\varphi \in \Phi_d$ has a hole effect, then $C(\bm{h}) = \varphi\left(\sqrt{\bm{h}^\top{\bf A}\bm{h}}\right)$ is a valid stationary covariance function, for any positive semidefinite matrix ${\bf A}$. This is one of the most utilized strategies to introduce anisotropy from an initial isotropic model, known as geometric (if $|{\bf A}|>0$, with $|\cdot|$ denoting the determinant of a square matrix) or zonal (if $| {\bf A} |=0$) anisotropy. Thus, hole effects and geometric/zonal anisotropies can coexist in a single family. However, this construction is overly rigid because the hole effect is constrained to occur in (almost) all directions with the same sharpness; of course, depending on the direction, the hole effect is attained at different ranges. 

\item Constructions of the form $C(\bm{h}) = \varphi_1(\|\bm{h}\|) \varphi_2(|h_i|)$, with $\varphi_1\in\Phi_d$, $\varphi_2\in\Phi_1$ and $h_i$ being the $i$th element of $\bm{h}$, can exhibit hole effects in directions that are close to the $i$-th axis, provided that $\varphi_2$ has a hole effect, see for instance \cite{LeBlevec}. This approach also produces a pattern that is quite rigid, where the interval of negative values in all directions exhibiting a hole effect (primarily, in orientations approximately parallel to the $i$-th axis) has a similar length regardless of the direction considered. 
\end{itemize}

\begin{figure}
    \centering
    \includegraphics[scale=0.27]{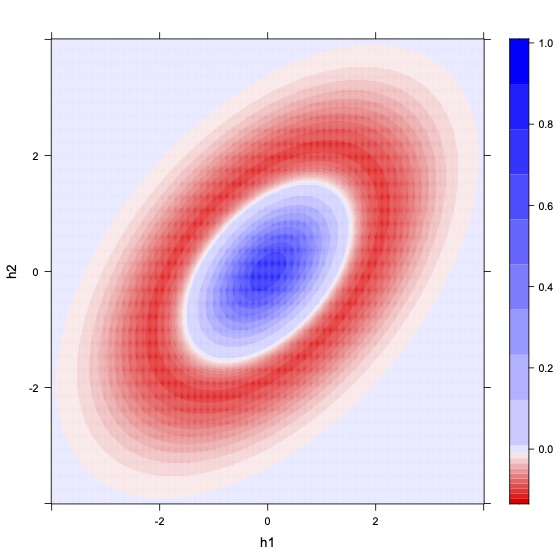}    \includegraphics[scale=0.27]{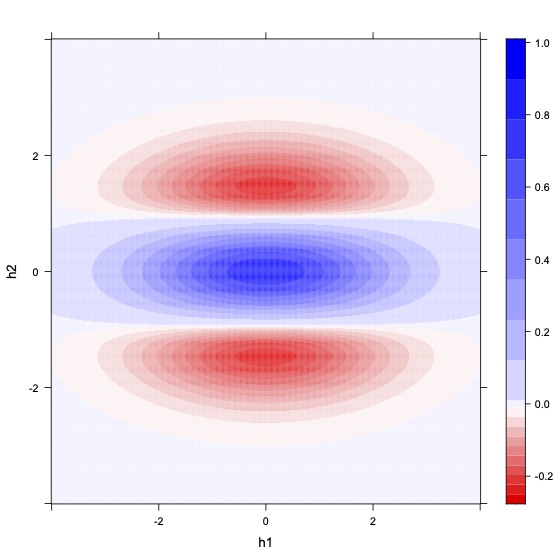}
\includegraphics[scale=0.27]{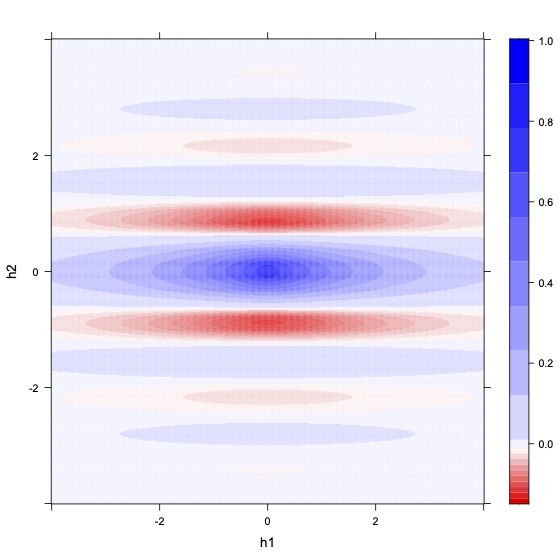}
    \caption{Basic constructions in dimension $d=2$. (Left) $2\exp(-0.8 \, \bm{h}^\top {\bf A} \bm{h}) -\exp(-0.4 \, \bm{h}^\top {\bf A} \bm{h})$, with ${\bf A} = {\scriptsize \begin{pmatrix}
        1&-0.5\\ -0.5&1
    \end{pmatrix}}$. (Middle) $\exp(-0.2 \, \|\bm{h}\|^2) \left[ 3.41 \exp(-0.8 \, h_2^2) -2.41 \exp(-0.4 \, h_2^2) \right]$. (Right) $\exp(-0.5 \, \|\bm{h}\|) \sin(5 |h_2|)/(5|h_2|)$. The positive semidefiniteness of the first two models, where differences of covariance functions are involved, is a consequence of Theorem 1(i) in \cite{ma2005linear}.}
    \label{fig:basic}
\end{figure}

Figure \ref{fig:basic} displays examples of these basic constructions, where the aforementioned structures can be visualized. This manuscript investigates other constructions that allow for complex combinations of these features.

\section{General Results}
\label{sec:results}

\subsection{Difference Between Geometrically Anisotropic Models}

In this section, we will examine the conditions under which the difference between two geometrically anisotropic covariance functions remains positive semidefinite.

\begin{proposition}
\label{prop1}
Let $\varphi$ be a member of the class $\Phi_d$ possessing a $d$-radial spectral density $f_d$.  Consider scalars $b_1, b_2 \geq 0$ and symmetric positive definite matrices ${\bf A}_1$ and ${\bf A}_2$. Thus, 
\begin{equation}
\label{covariance}
   \mathcal{T}^{(1)}_{{\bf A}_1,{\bf A}_2,b_1,b_2}[\varphi](\bm{h}) = 
      b_1 \, \varphi\left(\sqrt{\bm{h}^\top {\bf A}_1 \bm{h}}\right) - b_2  \, \varphi\left(\sqrt{\bm{h}^\top {\bf A}_2 \bm{h}}\right), \qquad \bm{h}\in\mathbb{R}^d,
\end{equation}
is a stationary covariance function in $\mathbb{R}^d$ if and only if
\begin{equation}
    \label{condition}
    b_1  \geq b_2 \, \frac{|{\bf A}_1|^{1/2}}{|{\bf A}_2|^{1/2}} \sup_{\bm{\omega}\in\mathbb{R}^d} \frac{f_d\left(\sqrt{\bm{\omega}^\top {\bf A}_2^{-1} \bm{\omega}}\right)}{ f_d\left(\sqrt{\bm{\omega}^\top {\bf A}_1^{-1} \bm{\omega}} \right)}.
\end{equation}
\end{proposition}

\begin{prueba}
Based on Bochner's theorem, one must show that the inverse Fourier transform of \eqref{covariance}, which is positively proportional to
\begin{equation}
\label{fourier}
  b_1 \int_{\mathbb{R}^d} 
 \exp\left(-\imath \bm{\omega}^\top \bm{h}\right) \varphi\left(\sqrt{\bm{h}^\top {\bf A}_1 \bm{h}}\right) \text{d}\bm{h} 
 - b_2 \int_{\mathbb{R}^d} 
 \exp\left(-\imath \bm{\omega}^\top \bm{h}\right) \varphi\left(\sqrt{\bm{h}^\top {\bf A}_2 \bm{h}}\right) \text{d}\bm{h},
\end{equation}
is nonnegative for every $\bm{\omega} \in \mathbb{R}^d$.
A change of variable allows writing \eqref{fourier} in the following format
\begin{equation*}
\frac{b_1}{|{\bf A}_1|^{1/2}} \int_{\mathbb{R}^d} \exp\left(-\imath  \left[{\bf A}_1^{-1/2} \bm{\omega}\right]^\top   \bm{v} \right)  
\varphi\left( \sqrt{\bm{v}^\top  \bm{v}}\right) \text{d}\bm{v} - 
\frac{b_2}{|{\bf A}_2|^{1/2}} \int_{\mathbb{R}^d} \exp\left(-\imath  \left[{\bf A}_2^{-1/2} \bm{\omega}\right]^\top   \bm{v} \right)  \varphi\left( \sqrt{\bm{v}^\top \bm{v}} \right) \text{d}\bm{v}.
\end{equation*}
Thus, up to a positive factor, \eqref{fourier} can be written as
\begin{equation}
\label{spectral_density}
\frac{b_1}{|{\bf A}_1|^{1/2}} f_d\left(\sqrt{\bm{\omega}^\top {\bf A}_1^{-1} \bm{\omega}}\right)  -  \frac{b_2}{|{\bf A}_2|^{1/2}}  f_d\left(\sqrt{\bm{\omega}^\top {\bf A}_2^{-1} \bm{\omega}}\right).
\end{equation}
The proof is completed by noting that \eqref{spectral_density} is nonnegative, for all $\bm{\omega} \in \mathbb{R}^d$, if and only if \eqref{condition} holds.     

\hfill  $\qed$
\end{prueba}

The term with a negative sign in \eqref{covariance} is the one that induces the hole effect, so matrix ${\bf A}_2$ is essential to characterize the predominant directions of the hole effect.

When the spectral density is radial and nonincreasing, the previous proposition can be simplified. Before stating the next result, we introduce the notation ${\bf A}_1 \succeq {\bf A}_2$, which indicates that ${\bf A}_1 - {\bf A}_2$ is a positive semidefinite matrix.

\begin{cor}
\label{cor1}
Let $\varphi$ be a member of the class $\Phi_d$ having a nonincreasing $d$-radial spectral density $f_d$. Let ${\bf A}_1$ and ${\bf A}_2$ be positive definite matrices such that ${\bf A}_1 \succeq {\bf A}_2$, and $b_1,b_2\geq 0$. 
 Thus, \eqref{covariance} is a stationary covariance function in $\mathbb{R}^d$ if and only if 
\begin{equation}
\label{simplified_condition}
b_1 \geq b_2 \,\frac{|{\bf A}_1|^{1/2}}{|{\bf A}_2|^{1/2}}.
\end{equation}
\end{cor}

\begin{prueba}
Condition ${\bf A}_1 \succeq {\bf A}_2$ is equivalent to ${\bf A}_2^{-1} \succeq {\bf A}_1^{-1}$. Thus, $\bm{\omega}^\top {\bf A}_2^{-1} \bm{\omega} \geq   \bm{\omega}^\top {\bf A}_1^{-1}   \bm{\omega}$ for all $\bm{\omega} \in \mathbb{R}^d$. Since $f_d$ is nonincreasing, 
$$f_d\left(\sqrt{\bm{\omega}^\top {\bf A}_2^{-1}   \bm{\omega}}\right)  \leq f_d\left(\sqrt{\bm{\omega}^\top {\bf A}_1^{-1}   \bm{\omega}}\right), \qquad \bm{\omega}\in\mathbb{R}^d.$$ 
Consequently, the supremum in the right hand side of \eqref{condition} is identically equal to one (attained for $\bm{\omega} = \bm{0}$).
$\hfill \qed$
\end{prueba}

\begin{remark}
\label{montee}
A sufficient condition for the $d$-radial spectral density $f_d$ to be nonincreasing is that $\varphi$ belongs to $\Phi_{d+2}$ and possesses a $(d+2)$-radial spectral density $f_{d+2}$. Indeed, in such a case, $\varphi$ is the Hankel transform of order $(d-2)/2$ of $f_d$, as per \eqref{fourier_inversion2}, and also the Hankel transform of order $d/2$ of $f_{d+2}$. This entails that $f_d$ is the mont\'ee of order $2$ of $f_{d+2}$ \citep[formula I.4.8]{matheron1965}:
\begin{equation}
f_d(u) = {2\pi} \int_u^{\infty} v f_{d+2}(v) \text{d}v, \quad u \geq 0,
\end{equation}
which is a nonincreasing function of $u$ insofar as $f_{d+2}$ is nonnegative.
\end{remark}

The conditions in the previous corollary can be stated in terms of the eigenvalues of ${\bf A}_1$ and ${\bf A}_2$. Let us denote by $\lambda_{j}({\bf A}_i)$, $\lambda_{\min}({\bf A}_i)$ and $\lambda_{\max}({\bf A}_i)$, the $j$-th, minimum and maximum eigenvalues of matrix ${\bf A}_i$, respectively, for $i=1,2$ and $j=1,\ldots, d$.

\begin{cor}
\label{cor2}
Let $\varphi$ be a member of the class $\Phi_d$ having a nonincreasing $d$-radial spectral density. Let ${\bf A}_1$ and ${\bf A}_2$ be positive definite matrices such that $\lambda_{\min}({\bf A}_1)  \geq  \lambda_{\max}({\bf A}_2)$, and $b_1,b_2\geq 0$.
 Thus, \eqref{covariance} is a stationary covariance function in $\mathbb{R}^d$ if and only if 
\begin{equation}
\label{condition_eigenvalues}
b_1 \geq b_2 \left( \prod_{j=1}^d  \frac{\lambda_{j}({\bf A}_1)}{\lambda_{j}({\bf A}_2)} \right)^{1/2}.
\end{equation}
\end{cor}

\begin{remark}
\label{remark1}
When ${\bf A}_i = a_i {\bf I}_{d}$, for $i=1,2$, with $a_{1} \geq a_2$ and ${\bf I}_d$ being the $d\times d$ identity matrix, 
 \eqref{covariance} reduces to the isotropic model
 \begin{equation}
     \label{nested}
     h \mapsto b_1\, \varphi(\sqrt{a_1}h) - b_2\, \varphi(\sqrt{a_2}h), 
 \end{equation}
 with $h= \|\bm{h}\|\geq 0$, 
and the respective validity condition \eqref{condition_eigenvalues} simplifies into 
\begin{equation}
\label{isotropic_condition}
    b_1 \geq b_2 \left(\frac{a_1}{a_2}\right)^{d/2}.
\end{equation}
Our results align with prior literature concerning this topic in the purely isotropic case. Specifically, we recover Theorem 1(ii) in \cite{ma2005linear}, and generalize Theorem 3.1 in \cite{faouzi2020} and Corollaries 3-12 in \cite{posa}. The results of this section can therefore be seen as an anisotropic extension of previous literature related to the difference between isotropic covariance models (or nested models) and the so-called Zastavnyi operators. 
\end{remark}

\subsection{Construction Based on Shifted Isotropic Models}

We propose here an alternative approach for constructing anisotropic covariance functions that exhibit negative values in specific orientations. We start with an isotropic model of the form \eqref{nested}. Therefore, it becomes crucial to satisfy both condition \eqref{isotropic_condition} and the requirement of having a nonincreasing $d$-radial spectral density for $\varphi$ to ensure that we start with an admissible covariance model. Then, we incorporate a shift in a determined direction to produce an anisotropic structure.

\begin{proposition}
\label{prop2}
Let $\varphi \in \Phi_d$ possessing a nonincreasing $d$-radial spectral density and consider constants $a_1,a_2 > 0$ and $b_1,b_2 \geq 0$ such that \eqref{isotropic_condition} holds. Thus, for all $\bm{\eta}\in\mathbb{R}^d$, the mapping 
\begin{equation}
\label{construction2} \mathcal{T}^{(2)}_{a_1,a_2,b_1,b_2,\bm{\eta}}[\varphi](\bm{h}) =  b_1 \,   \varphi\left(\sqrt{a_1} \|\bm{h}\|\right) - \frac{b_2}{2} \big[ \varphi(\sqrt{a_2} \|\bm{h}-\bm{\eta}\|) + \varphi(\sqrt{a_2} \|\bm{h}+\bm{\eta}\|) \big], \qquad \bm{h}\in\mathbb{R}^d, 
\end{equation}
is a stationary covariance function in $\mathbb{R}^d$.
\end{proposition}

\begin{prueba}
Let $f_{a_i,d}$ denote the $d$-radial spectral density of $\varphi(\sqrt{a_i} h)$, for $i=1,2$. Note that
\begin{equation*}
\begin{split}
   \frac{1}{(2\pi)^d} \int_{\mathbb{R}^d} 
 \exp\left(-\imath \bm{\omega}^\top \bm{h}\right) \varphi\left(\sqrt{a_2} \|\bm{h}-\bm{\eta}\|\right) \text{d}\bm{h} &= \frac{1}{(2\pi)^d}\int_{\mathbb{R}^d} 
 \exp\left(-\imath \bm{\omega}^\top \left[\bm{v}+ \bm{\eta}\right]\right) \varphi\left(\sqrt{a_2} \|\bm{v}\|\right) \text{d}\bm{v} \\
 &= \exp\left(-\imath \bm{\omega}^\top \bm{\eta}\right)f_{a_2,d}(\omega),
\end{split}    
\end{equation*}
 for all $\bm{\omega}\in\mathbb{R}^d$, with $\omega = \|\bm{\omega}\|$. Similarly,
 $$  \frac{1}{(2\pi)^d}\int_{\mathbb{R}^d} 
 \exp\left(-\imath \bm{\omega}^\top \bm{h}\right) \varphi\left(\sqrt{a_2} \|\bm{h}+\bm{\eta}\|\right) \text{d}\bm{h} = \exp\left(\imath \bm{\omega}^\top \bm{\eta}\right)f_{a_2,d}(\omega). $$
 Thus, the inverse Fourier transform of \eqref{construction2} can be written as
\begin{equation} 
\label{fourier2}
b_1 f_{a_1,d}(\omega) - \frac{b_2}{2} \left[   \exp\left(-\imath \bm{\omega}^\top \bm{\eta}\right)f_{a_2,d}(\omega) + \exp\left(\imath \bm{\omega}^\top \bm{\eta}\right)f_{a_2,d}(\omega)
 \right] =  b_1 f_{a_1,d}(\omega) - b_2   \cos\left(\bm{\omega}^\top \bm{\eta}\right) f_{a_2,d}(\omega) 
\end{equation}
for all $\bm{\omega}\in\mathbb{R}^d$. The right-hand side of \eqref{fourier2} is lower-bounded by $b_1 f_{a_1,d}(\omega) - b_2 f_{a_2,d}(\omega)$,  where the latter expression corresponds to the $d$-radial spectral density of \eqref{nested}. This quantity is non-negative because condition \eqref{isotropic_condition} is satisfied, i.e., \eqref{nested} is positive semidefinite. The proof is completed by invoking Bochner's theorem. 
\hfill $\qed$
\end{prueba}

The interest of the above proposition lies in the fact that all the isotropic constructions of the form \eqref{nested} can be adapted according to \eqref{construction2} to produce anisotropic models. When the separation vector $\bm{h}$ is close to $\pm \bm{\eta}$, the negative part of \eqref{construction2} becomes  predominant; thus, the hole effect is more significant in that direction. 

There are two limit cases of \eqref{construction2} worth noting. On the one hand, as the magnitude of $\bm{\eta}$ approaches infinity, \eqref{construction2} tends to $b_1 \, \varphi(\sqrt{a_1} h)$ (a rescaled version of the initial covariance model). On the other hand, when the magnitude of $\bm{\eta}$ approaches zero, the nested model \eqref{nested} is recovered. Thus, this construction can encompass purely isotropic models, both with and without hole effect, as special cases.

\subsection{Models with Derivative Information}

Our focus now turns to the study of anisotropic models whose construction incorporates directional derivatives of an isotropic random field. In contrast to previous strategies, this approach requires a covariance function twice differentiable at the origin as one of the initial ingredients, and no monotonicity conditions are required for the $d$-radial spectral density.

\begin{proposition}
\label{prop3}
    Let $\varphi_1,\varphi_2 \in \Phi_d$, with $\varphi_2$ being twice differentiable at the origin, and $\bm{u}$ be a unit vector in $\mathbb{R}^d$. Consider constants $a_1,a_2 > 0$ and $b_1,b_2 \geq 0$. Thus, the mapping
    \begin{equation}
        \label{construction4}  \mathcal{T}^{(3)}_{a_1,a_2,b_1,b_2,\bm{u}}[\varphi_1,\varphi_2](\bm{h}) = b_1 \varphi_1(\sqrt{a_1}\|\bm{h}\|) - b_2 \left[ \cos^2(\theta(\bm{h},\bm{u}))\varphi_2''(\sqrt{a_2}\|\bm{h}\|) +  \sin^2(\theta(\bm{h},\bm{u}))\frac{\varphi_2'(\sqrt{a_2}\|\bm{h}\|)}{\sqrt{a_2}\|\bm{h}\|} \right],
    \end{equation}
    where $\bm{h}\in\mathbb{R}^d$, 
    with $\theta(\bm{h},\bm{u})$ being the angle between $\bm{h}$ and $\bm{u}$,
    is a stationary covariance function in $\mathbb{R}^d$.
\end{proposition}

\begin{prueba}
    We provide a constructive proof. Let us consider two independent zero-mean random fields on $\mathbb{R}^d$, denoted as $Y_1$ and $Y_2$, which possess covariance functions $\varphi_1$ and $\varphi_2$ in $\Phi_d$, respectively.  
    Equation (5.29) in \cite{chiles2009geostatistics} establishes that
    $${\rm cov}\left[ \frac{\partial {Y}_2}{\partial \bm{u}}(\bm{x}), \frac{\partial {Y}_2}{\partial \bm{v}}(\bm{x}+
    \bm{h})  \right]  = - \frac{\left(\bm{h}^\top \bm{u}\right)^2}{\|\bm{h}\|^2} \left[ \varphi_2''(\|\bm{h}\|) - \frac{\varphi_2'(\|\bm{h}\|)}{\|\bm{h}\|} \right] - \left(\bm{u}^\top \bm{v}\right) \frac{\varphi_2'(\|\bm{h}\|)}{\|\bm{h}\|}, $$
    for all $\bm{x}, \bm{h}\in\mathbb{R}^d$ and any pair of unit vectors $\bm{u}$ and $\bm{v}$ in $\mathbb{R}^d$, provided that $\varphi_2$ is twice differentiable at the origin. 
    Thus, a direct calculation shows that the covariance function of the random field $\left\{(\partial {Y}_2/\partial \bm{u})( \bm{x}): \bm{x}\in\mathbb{R}^d\right\}$ is given by
$$ \bm{h} \mapsto -   \cos^2(\theta(\bm{h},\bm{u}))  \varphi_2''(\|\bm{h}\|) -\sin^2(\theta(\bm{h},\bm{u}))\frac{\varphi_2'(\|\bm{h}\|)}{\|\bm{h}\|}. $$
    Based on previous calculations, one concludes that a random field defined according to
    \begin{equation}
        \label{derivatives}
        Z(\bm{x}) = \sqrt{b_1} Y_1(\sqrt{a_1} \bm{x}) + \sqrt{\frac{b_2}{a_2}}  \frac{\partial {Y}_2}{\partial \bm{u}}(\sqrt{a_2} \bm{x}), \qquad \bm{x}\in\mathbb{R}^d,
    \end{equation}
    has a covariance function given by \eqref{construction4}, indicating that \eqref{construction4} is positive semidefinite.
    \hfill $\qed$
\end{prueba}

The rationale behind this approach is that the changes in sign of the directional derivative in \eqref{derivatives} can accentuate the transitions between large and small values of the random field $Z$ in a given direction; thus, marked hole effects in the orientation determined by $\bm{u}$ are expected. If $\bm{h}$ is approximately proportional to $\bm{u}$, the second-order derivative of $\varphi_2$ gains greater significance in \eqref{construction4}. Conversely, if $\bm{h}$ is approximately orthogonal to $\bm{u}$, the term involving the first-order derivative becomes more dominant.

The parameters involved in this formulation do not require any elaborate restriction, as the positive semidefiniteness is inherently ensured by construction. A special case of \eqref{construction4} arises when setting $b_1 = 0$, where the dominant component of the covariance structure is the term within brackets, representing the covariance function of the directional derivative of certain random field.

When the covariance functions of $Y_1$ and $Y_2$ are equal and given by $\varphi_1 = \varphi_2 := \varphi$, where $\varphi$ is a function in $\Phi_d$ that is twice differentiable at the origin, we can conveniently denote the expression \eqref{construction4} as $\mathcal{T}^{(3)}_{a_1,a_2,b_1,b_2,\bm{u}}[\varphi]$.

\begin{remark}
    \label{remark:construction3}
 It is noteworthy that, in Proposition \ref{prop3}, one can substitute $\varphi_1$ with a stationary covariance model, which need not be isotropic. The validity of this alternative model is guaranteed by following the same proof as before. This slight variation offers enhanced flexibility in spatial data modeling.
\end{remark}

\section{Explicit Parametric Families}
\label{parametric}

\subsection{Mat\'ern, Cauchy and Compactly-Supported Hypergeometric Models}

To provide concrete models derived from the findings presented in the previous section, we will now introduce three commonly used parametric families of covariance functions: the Mat\'ern, Cauchy and Gauss hypergeometric families.
\begin{enumerate}
    \item The Mat\'ern family of covariance functions is given by \citep{stein-book}
\begin{equation}
    \label{matern}
    \mathcal{M}_{\nu}(t) = \frac{2^{1-\nu}}{\Gamma(\nu)} t^\nu \mathcal{K}_{\nu}(t), \qquad t\geq 0,
\end{equation}
where 
$\mathcal{K}_\nu$ is the modified Bessel function of the second kind, with $\nu>0$ being a shape parameter \citep{Olver}. The $d$-radial spectral density associated with this model, viewed as a function of $\omega = \|\bm{\omega}\|$, is given by 
$$f_d^{\mathcal{M}}(\omega) = \frac{\Gamma(\nu+d/2)}{\Gamma(\nu) \pi^{d/2}} \frac{1}{(1+\omega^2)^{\nu + d/2}}, \qquad \omega\geq 0. $$
\item The Cauchy family of covariance functions is given by (see, e.g., \citealp{chiles2009geostatistics})
\begin{equation}
    \label{matern}
    \mathcal{C}_{\delta}(t) = (t^2+1)^{-\delta}, \qquad t\geq 0,
\end{equation}
with $\delta>0$ being a shape parameter. When $\delta > (d-1)/4$, its $d$-radial spectral density adopts the explicit form \citep{lim2009gaussian}
$$ f_d^{\mathcal{C}}(\omega) = \frac{2^{1-d/2-\delta}}{\Gamma(\delta) \pi^{d/2}}  \frac{\mathcal{K}_{d/2-\delta}(\omega)}{\omega^{d/2-\delta}}, \qquad \omega\geq 0. $$
\item The Gauss hypergeometric family of covariance functions is given by \citep{emery2021gauss}
\begin{equation}
    \label{hygeo}
    \mathcal{H}_{\alpha,\beta,\gamma}(t) = (1-t^2)_+^{\beta-\alpha+\gamma-d/2-1} {}_2F_1(\beta-\alpha,\gamma-\alpha;\beta-\alpha+\gamma-d/2;(1-t^2)_+), \qquad t\geq 0,
\end{equation}
with ${}_2F_1$ denoting the Gauss hypergeometric function \citep{Olver},  $(\cdot)_+$ denoting the positive part and $\alpha, \beta, \gamma$ being shape parameters such that $2\alpha>d$, $2(\beta-\alpha)(\gamma-\alpha) \geq \alpha$ and $2(\beta+\gamma) \geq 6\alpha+1$. Its $d$-radial spectral density is 
$$ f_d^{\mathcal{H}}(\omega) = \kappa(\alpha;\beta,\gamma) {}_1F_2(\alpha;\beta,\gamma;-\omega^2/2), \qquad \omega\geq 0, $$
with $\kappa(\alpha;\beta,\gamma)$ a positive factor and ${}_1F_2$ a generalized hypergeometric function \citep{Olver}. This model encompasses the Euclid's hat (spherical), cubic, generalized Wendland and Askey covariances as particular cases. 
\end{enumerate}
Both $\mathcal{M}_\nu$ and $\mathcal{C}_{\delta}$ belong to the class $\Phi_d$, for all $d\geq 1$, and both $f_d^{\mathcal{M}}$ and $f_d^{\mathcal{C}}$ are decreasing functions. As for $\mathcal{H}_{\alpha,\beta,\gamma}$, it belongs to $\Phi_{d+2}$ if $2\alpha>d+2$, $2(\beta-\alpha)(\gamma-\alpha) \geq \alpha$ and $2(\beta+\gamma) \geq 6\alpha+1$, in which case $f_d^{\mathcal{H}}$ is a nonincreasing function (recall Remark \ref{montee}). Thus, these three models are in the range of applicability of Propositions \ref{cor1} and \ref{prop2}. 

While the Cauchy model is infinitely differentiable at the origin \citep{chiles2009geostatistics} and so is the Gauss hypergeometric model if $2\alpha > d+2$ \citep{emery2021gauss}, the Mat\'ern model is twice differentiable at the origin if and only if $\nu>1$ \citep{stein-book} and, in this case, Proposition \ref{prop3} can be applied. 

In summary, we have the following corollaries. 

\begin{cor}
    \label{cor1_matern}
    Consider two positive definite matrices ${\bf A}_1$ and ${\bf A}_2$ such that ${\bf A}_1 \succeq {\bf A}_2$, and scalars $b_1,b_2\geq 0$. Thus, $\mathcal{T}^{(1)}_{{\bf A}_1,{\bf A}_2,b_1,b_2}[\mathcal{M}_\nu]$, $\mathcal{T}^{(1)}_{{\bf A}_1,{\bf A}_2,b_1,b_2}[\mathcal{C}_\delta]$ and $\mathcal{T}^{(1)}_{{\bf A}_1,{\bf A}_2,b_1,b_2}[\mathcal{H}_{\alpha,\beta,\gamma}]$,  with $\nu>0$, $\delta > (d-1)/4$, $2\alpha>d+2$, $2(\beta-\alpha)(\gamma-\alpha) \geq \alpha$ and $2(\beta+\gamma) \geq 6\alpha+1$, are stationary covariance functions in $\mathbb{R}^d$ 
    if and only if condition \eqref{simplified_condition} holds.
\end{cor}

\begin{cor}
    \label{cor2_matern}
Let $a_1,a_2 > 0$ and $b_1,b_2 \geq 0$ be constants satisfying condition \eqref{isotropic_condition} and $\bm{\eta} \in \mathbb{R}^d$. Thus, $\mathcal{T}^{(2)}_{a_1,a_2,b_1,b_2,\bm{\eta}}[\mathcal{M}_\nu]$, $\mathcal{T}^{(2)}_{a_1,a_2,b_1,b_2,\bm{\eta}}[\mathcal{C}_\delta]$ and $\mathcal{T}^{(2)}_{a_1,a_2,b_1,b_2,\bm{\eta}}[\mathcal{H}_{\alpha,\beta,\gamma}]$, with $\nu>0$, $\delta > (d-1)/4$, $2\alpha>d+2$, $2(\beta-\alpha)(\gamma-\alpha) \geq \alpha$ and $2(\beta+\gamma) \geq 6\alpha+1$, are stationary covariance functions in $\mathbb{R}^d$.
\end{cor}

\begin{cor}
    \label{cor3_matern}
Consider constants $a_1,a_2 > 0$ and $b_1,b_2 \geq 0$, and a unit vector $\bm{u} \in \mathbb{R}^d$. Thus, $\mathcal{T}^{(3)}_{a_1,a_2,b_1,b_2,\bm{u}}[\mathcal{M}_\nu]$ with $\nu>1$,  $\mathcal{T}^{(3)}_{a_1,a_2,b_1,b_2,\bm{u}}[\mathcal{C}_\delta]$ and $\mathcal{T}^{(3)}_{a_1,a_2,b_1,b_2,\bm{u}}[\mathcal{H}_{\alpha,\beta,\gamma}]$ with  $2\alpha>d+2$, $2(\beta-\alpha)(\gamma-\alpha) \geq \alpha$ and $2(\beta+\gamma) \geq 6\alpha+1$,
are stationary covariance functions in $\mathbb{R}^d$.
\end{cor}

In order to exhibit the versatility of the proposed models, we provide visual illustrations in dimension $d=2$. These illustrations show the various shapes that can be achieved. We consider the following scenarios:
\begin{itemize}
    \item[{\bf I.}] The models in Corollary \ref{cor1_matern}, with ${\bf A}_1 = {\bf I}_{2}$  and ${\bf A}_2 = {\bf P} \, \text{diag}(\mu_1,\mu_2) \, {\bf P}^\top$, with $\mu_1,\mu_2 > 0$ and 
    $${\bf P} = \begin{bmatrix}
\cos(\pi/4) & -\sin(\pi/4)\\
\sin(\pi/4) & \cos(\pi/4)
\end{bmatrix}$$
    being a rotation matrix. The conditions of Corollary \ref{cor1_matern} are satisfied if and only if $\max(\mu_1,\mu_2) \leq 1$ and $ b_1  \sqrt{\mu_1 \mu_2} \geq b_2$. Thus, we fix $b_1 = 2.5$, $b_2 = 1$, $\mu_1 = 0.2$ and $\mu_2=0.8$.
\item[{\bf II.}] The models in Corollary \ref{cor2_matern}, with $b_1 = 2$, $b_2=1$, $a_1 = 0.8$ and $a_2=0.4$, with a shift vector given by $\bm{\eta} = [1,1]^\top$.  

\item[{\bf III.}] The models in Corollary \ref{cor3_matern}, with 
 $b_1=1$, $b_2=2$, $a_1=1$ and $a_2 = 0.5$, and the unit vector $\bm{u} = [1/\sqrt{2},1/\sqrt{2}]^\top$.  
\end{itemize}

\begin{figure}
\centering
\includegraphics[scale=0.27]{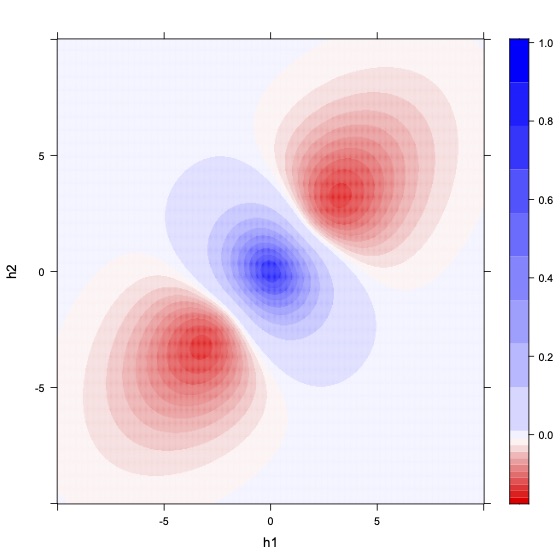}
\includegraphics[scale=0.27]{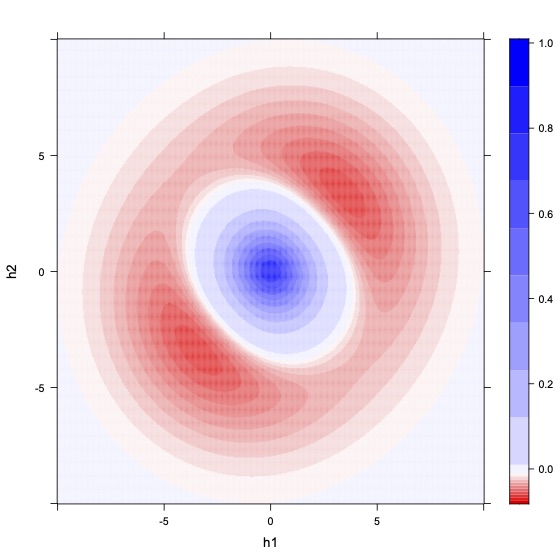}
\includegraphics[scale=0.27]{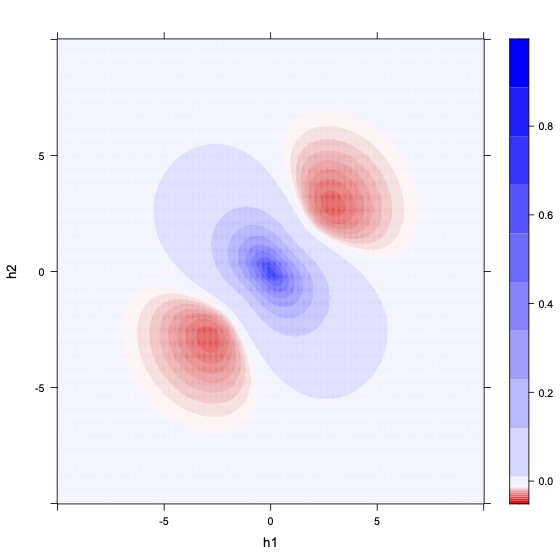}

\includegraphics[scale=0.27]{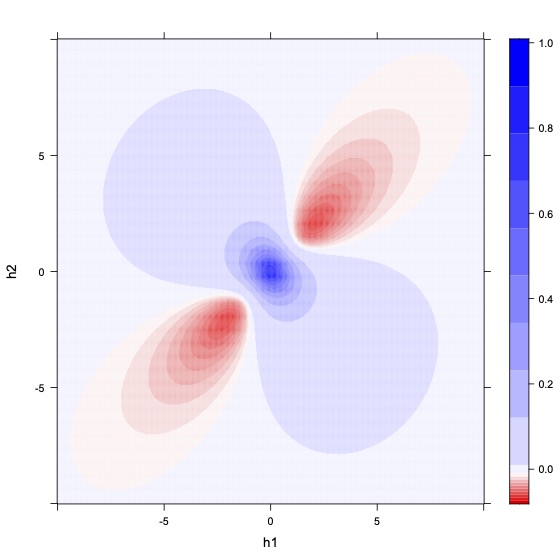}
\includegraphics[scale=0.27]{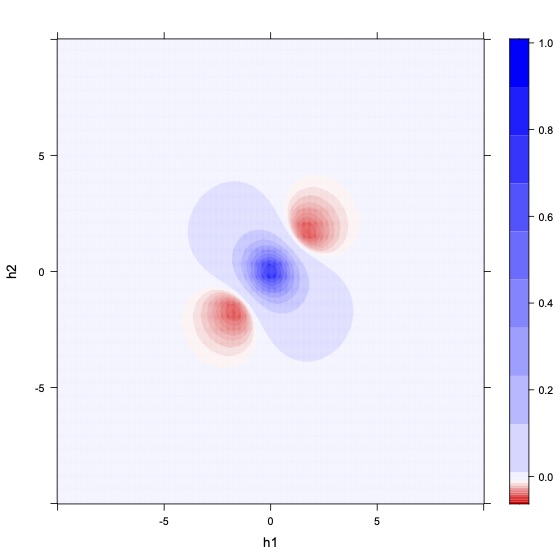}
\includegraphics[scale=0.27]{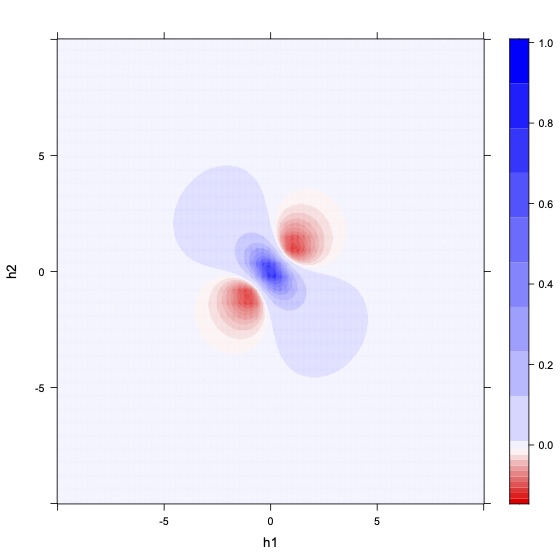}

\includegraphics[scale=0.27]{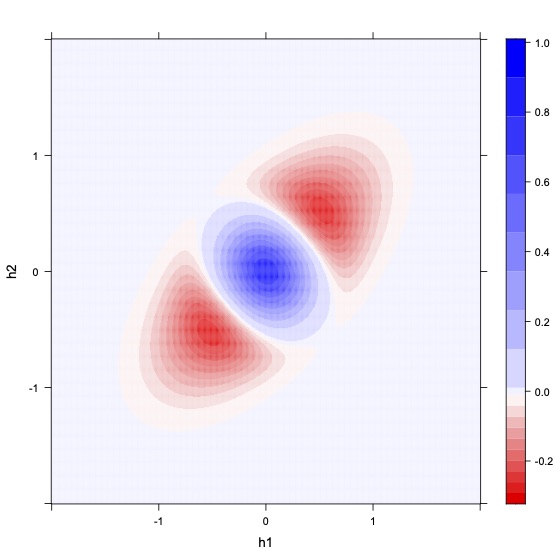}
\includegraphics[scale=0.27]{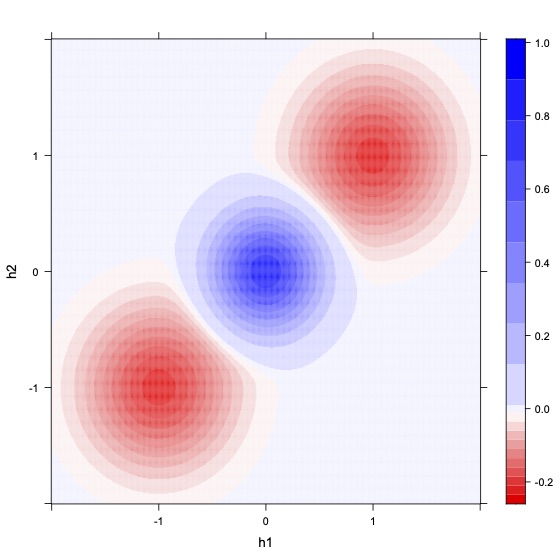}
\includegraphics[scale=0.27]{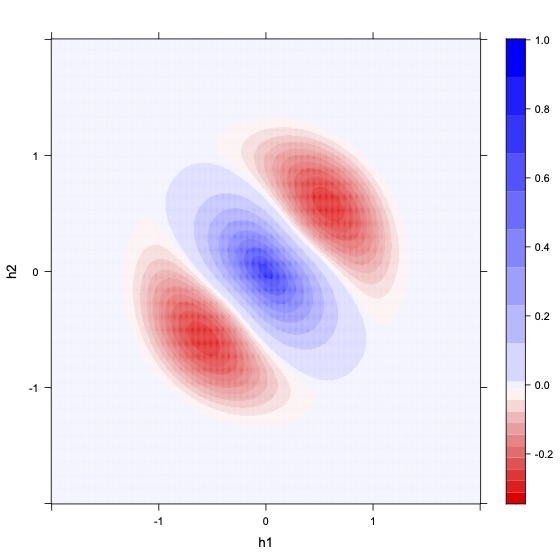}
\caption{Different combinations of anisotropies and hole-effects for the transformed Mat\'ern (top), the transformed Cauchy {(middle) and the transformed Gauss hypergeometric (bottom)} models. From left to right we consider the transformations introduced in Corollaries \ref{cor1_matern}-\ref{cor3_matern}, respectively. The values of the parameters have been described in scenarios {\bf I}-{\bf III}.}
\label{contours}
\end{figure}

Figure \ref{contours} shows the contour plots of the Mat\'ern model with $\nu=1.5$, the Cauchy model with $\delta = 1$ {and the Gauss hypergeometric model with $\alpha = 3, \beta = 7/2$ and $\gamma = 6$}, after the application of the transformations described in Corollaries \ref{cor1_matern}-\ref{cor3_matern} under scenarios {\bf I}-{\bf III}, respectively, together with a normalization in order to obtain correlation functions. To improve the visualization of each individual model, we have chosen specific ranges for plotting. We consider $\bm{h}=[h_1,h_2]^\top \in [-10,10]^2$ for the first two models, and $\bm{h}=[h_1,h_2]^\top \in [-2,2]^2$ for the last model. All the covariance functions have been designed to present a hole effect around the northeast direction.

\subsection{Cardinal Sine Model}

Our focus now turns to the cardinal sine (or wave) covariance function, defined through
\begin{equation}
    \label{wave}
    \mathcal{W}(t) = \frac{\sin(t)}{t}, \qquad t> 0,
\end{equation}
and $\mathcal{W}(0) = 1$. This model is a member of $\Phi_d$, for $d\leq 3$. When $d=3$, this model does not possess a spectral density. However, for $d\leq 2$, one has \citep{arroyo2021}
\begin{equation}
    \label{spectral_wave}
    f_d^{\mathcal{W}}(\omega) = \frac{1}{2 \pi^{(d-1)/2} \Gamma((3+d)/2)} (1-\omega^2)_+^{(1-d)/2}, \qquad \omega\geq 0.
\end{equation}
 In particular, when $d=2$ and $0 \leq  \omega < 1$, \eqref{spectral_wave} is an increasing mapping.  
As a result, Propositions \ref{prop1} and \ref{prop2} are not applicable to this model.
The conditions of Proposition \ref{prop3}, on the other hand, can be readily verified for $d\leq 3$, leading to the subsequent corollary.

\begin{cor}
\label{cor_wave}
Let $d\leq 3$. Consider constants $a_1,a_2 > 0$ and $b_1,b_2 \geq 0$, and a unit vector $\bm{u} \in \mathbb{R}^d$. Thus, $\mathcal{T}^{(3)}_{a_1,a_2,b_1,b_2,\bm{u}}[\mathcal{W}]$ is a stationary covariance function in $\mathbb{R}^d$.
\end{cor}

Recall that Proposition \ref{prop3} offers the flexibility to combine models from different parametric families. As an example, we can consider $\mathcal{T}^{(3)}_{a_1,a_2,b_1,b_2,\bm{u}}[\mathcal{M}_\nu,\mathcal{W}]$, which constitutes a valid stationary covariance model for dimensions $d\leq 3$. 

\begin{figure}
    \centering
\includegraphics[scale=0.27]{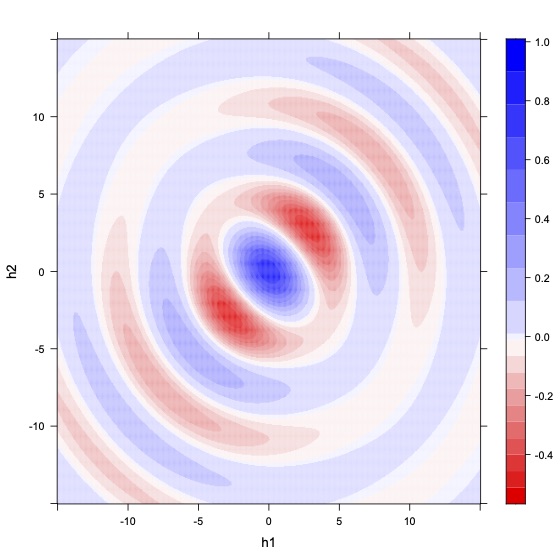} \includegraphics[scale=0.27]{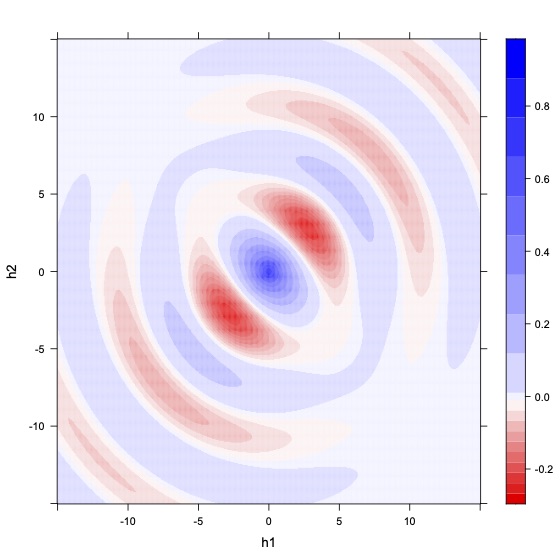}
    \caption{Models $\mathcal{T}^{(3)}_{a_1,a_2,b_1,b_2,\bm{u}}[\mathcal{W}]$ (Left) and $\mathcal{T}^{(3)}_{a_1,a_2,b_1,b_2,\bm{u}}[\mathcal{M}_{1/2},\mathcal{W}]$ (Right), with $d=2$,  
 $a_1=a_2=b_1=1$, $b_2=2$ and $\bm{u}=[1/\sqrt{2},1/\sqrt{2}]^\top$.}
    \label{contours_wave}
\end{figure}

Figure \ref{contours_wave} shows $\mathcal{T}^{(3)}_{a_1,a_2,b_1,b_2,\bm{u}}[\mathcal{W}]$ and $\mathcal{T}^{(3)}_{a_1,a_2,b_1,b_2,\bm{u}}[\mathcal{M}_{1/2},\mathcal{W}]$ in dimension $d=2$, with parameters $a_1=a_2=b_1=1$, $b_2=2$ and $\bm{u}=[1/\sqrt{2},1/\sqrt{2}]^\top$. While certain structural oscillations from the model \eqref{wave} persist, the proposed models exhibit a notably amplified hole effect in the $\bm{u}$ direction. Observe that $\mathcal{T}^{(3)}_{a_1,a_2,b_1,b_2,\bm{u}}[\mathcal{W}]$ exceeds the lower bound required for isotropic models in $\mathbb{R}^2$.

\section{Real Data Analysis}
\label{sec:data}

We consider a geophysical data set from a carbonate-rock aquifer located in Martin county, south Florida, and documented in \cite{parra2006, parra2009}. The data set consists of a P-wave impedance vertical section obtained by inverting cross-well reflection seismic measurements, at a vertical resolution of 0.61 m (2 feet) and a horizontal resolution of 3.05 m (10 feet), totaling 17,145 data. The P-wave impedance can be used to delineate the lateral heterogeneities of the aquifer, to assess the fluid paths, and to map petrophysical properties such as the rock porosity, which is a key variable to forecast water production \citep{parra2013, emery2013}.

\begin{figure}
    \centering
    \includegraphics[scale=0.41]{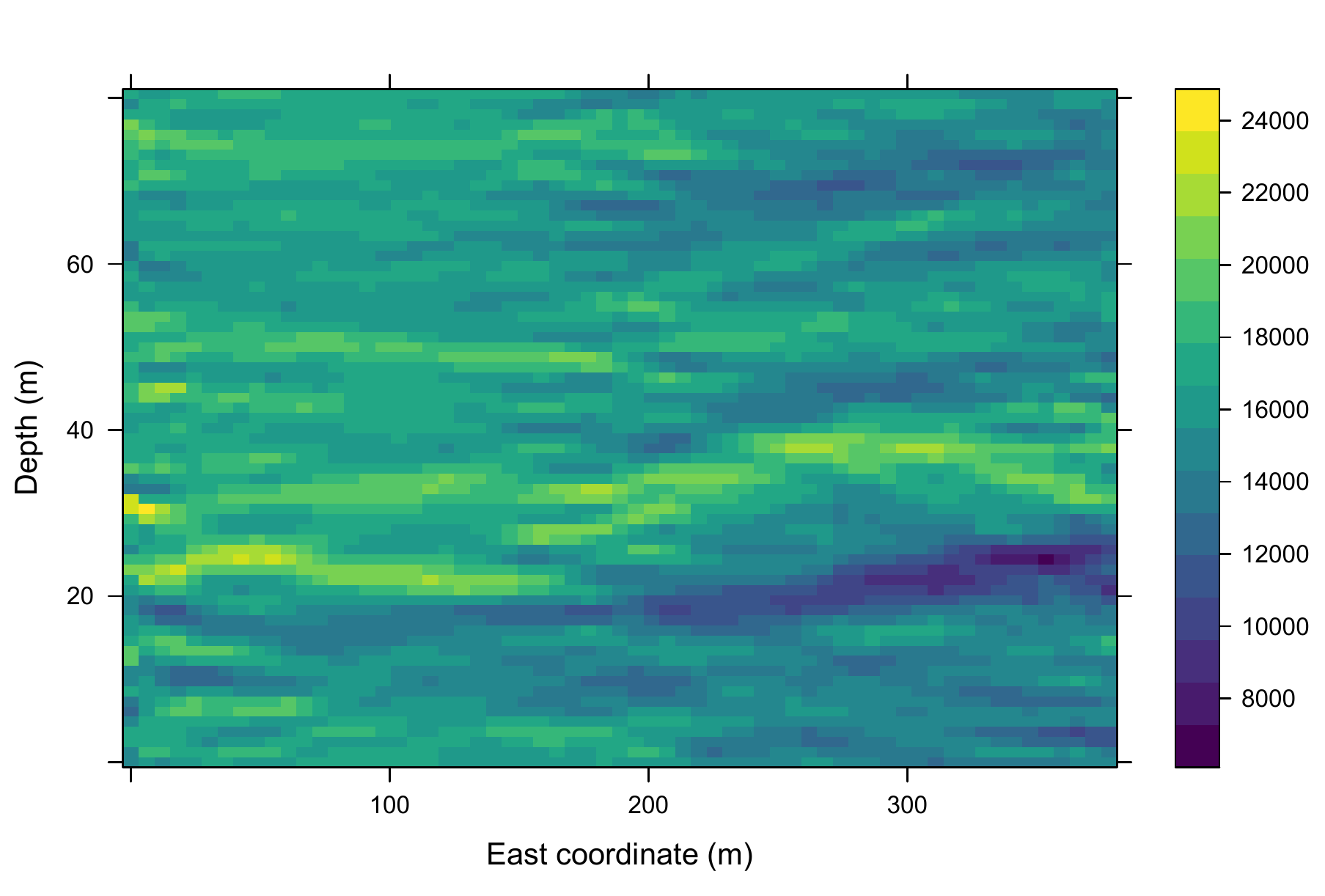} 
\includegraphics[scale=0.33]{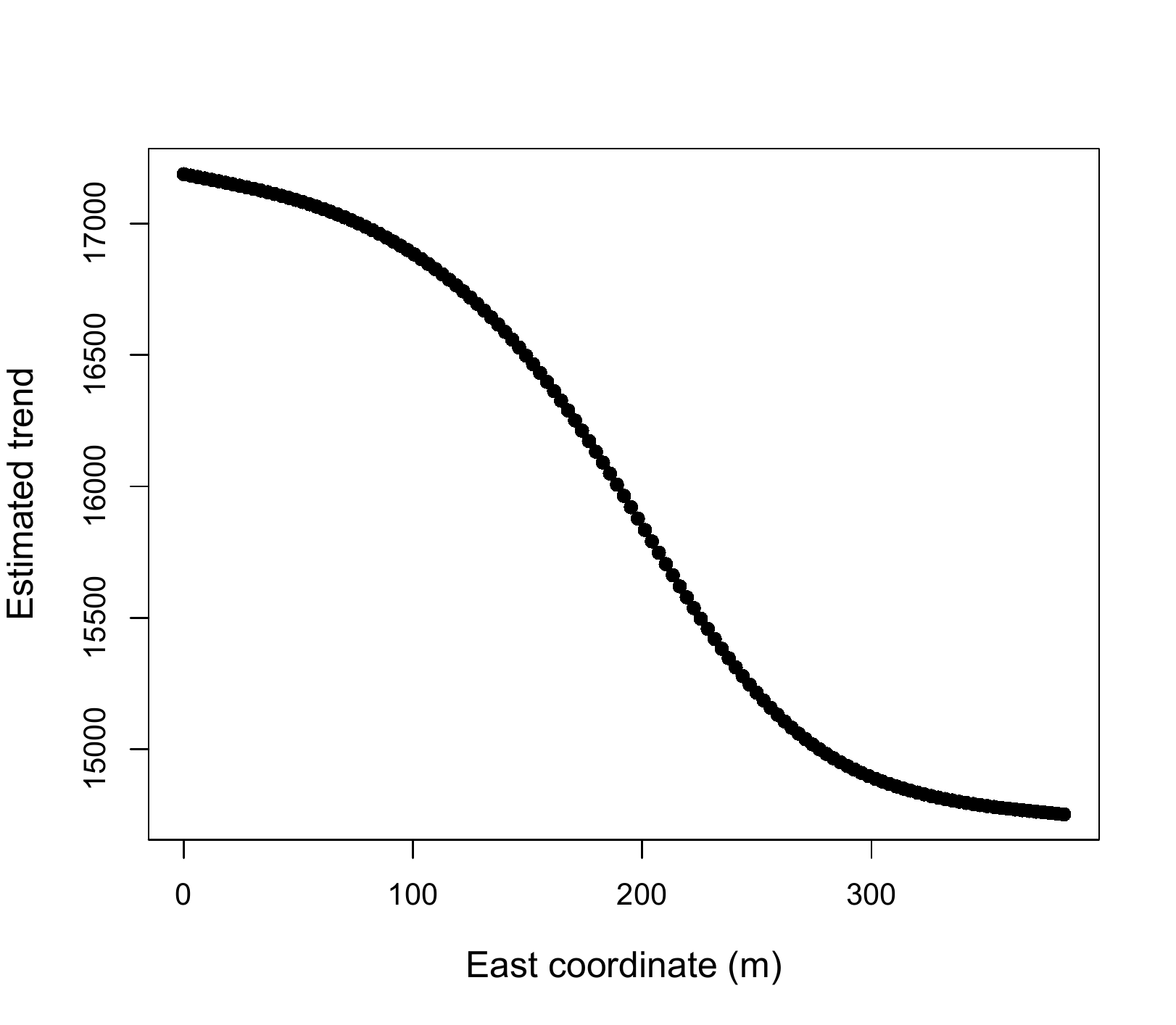}
        \includegraphics[scale=0.41]{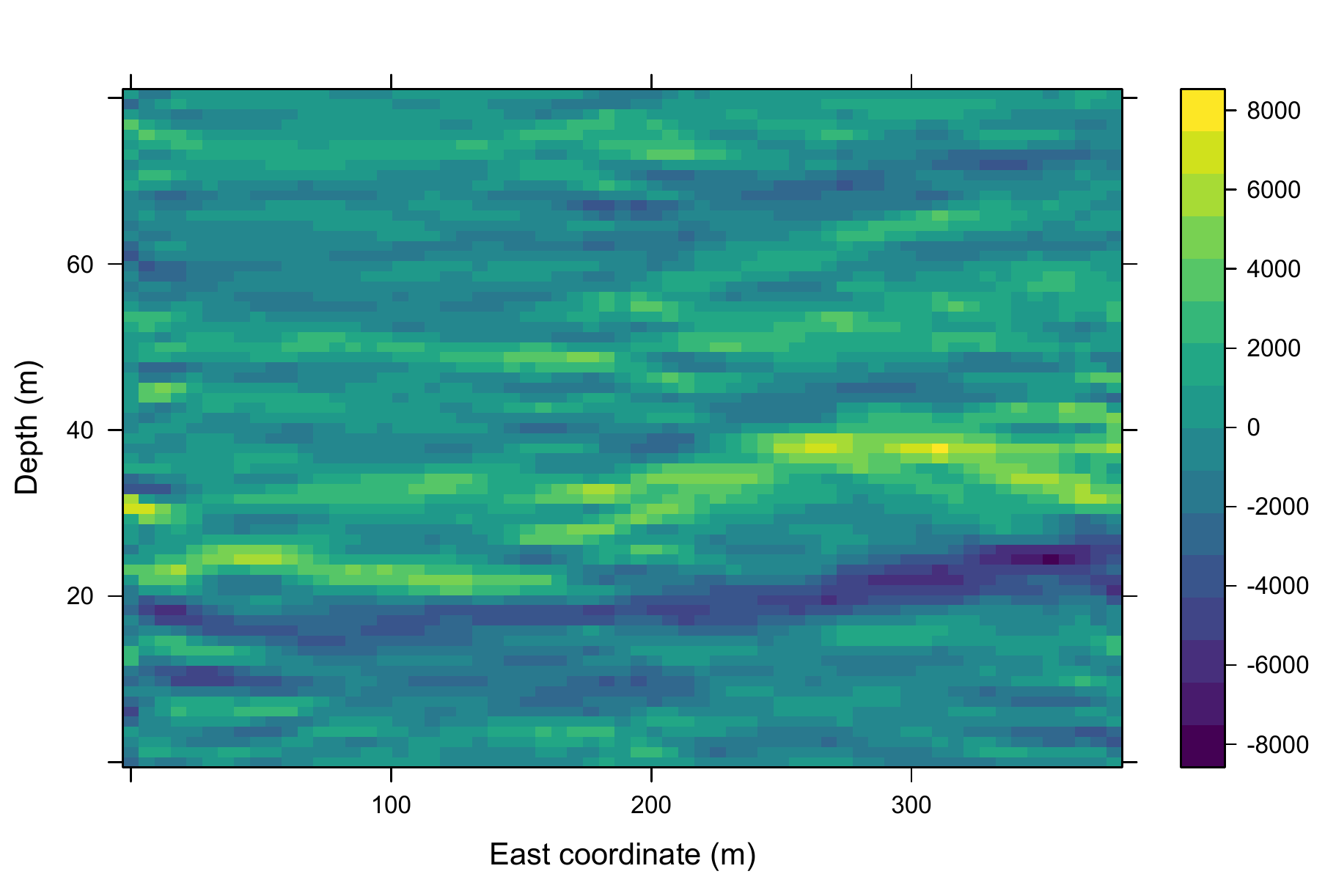} 
\includegraphics[scale=0.31]{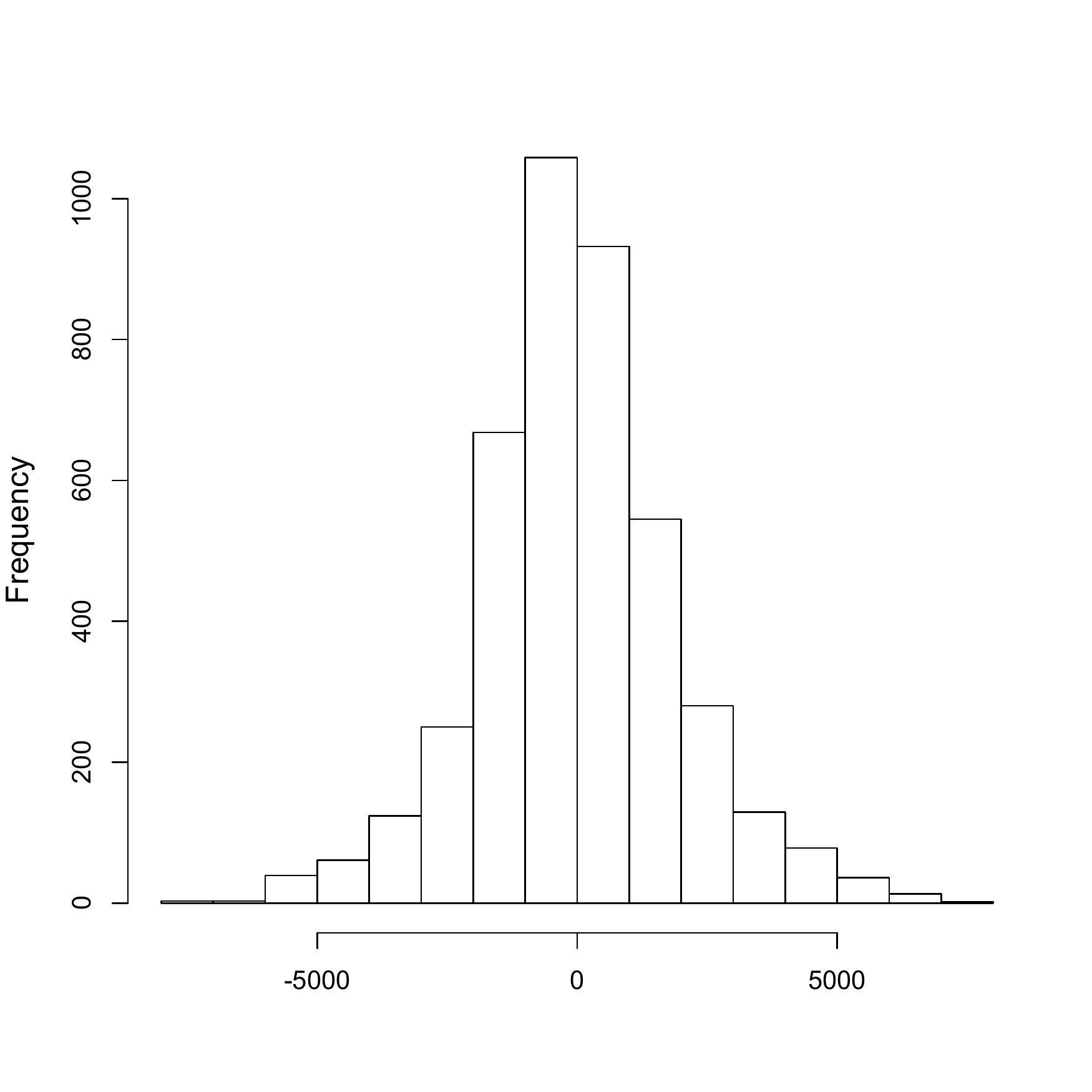}
    \caption{ From top left to bottom right: original data set of impedance, fitted trend in the east direction,
residuals and the corresponding histogram.    
}
    \label{fig:datos}
\end{figure}

To reduce the number of data, we employ for our analysis a spatial resolution of $20$ feet and $4$ feet in the horizontal and vertical coordinates, respectively, which leads to a set of 4352 impedance data. Also, to remove the trend in the east coordinate and improve the description of the data by a stationary random field model, we utilize a smoothing spline approach. 
The estimated trend exhibits a distinct pattern, gradually transitioning from high to low values as one moves from west to east.
In Figure \ref{fig:datos}, one can observe the original data, the trend that was fitted, the residuals, and the corresponding histogram. These residuals can be interpreted as the realization of a stationary zero-mean Gaussian random field. We randomly select and exclude $400$ observations of the dataset (approximately $10\%$ of the observations) for posterior validation purposes, while the remaining observations constitute the training set.

A significant hole effect is present in the vertical direction. This hole effect can be explained by the presence of major geological structures, corresponding to permeability barriers alternating vertically with high-porosity structures. The former are characterized by tight limestone and isolated vugs, whereas the latter are associated with interconnected matrix and vugs or with a combination of interconnected vugs surrounded by limestone \citep{parra2009}. Cyclic behaviors in the vertical covariances or variograms of rock properties are often observed in carbonate sequences and can be explained by periodic processes of deposition due to eustatic sea level oscillations or to tectonic activities \citep{chiles2009geostatistics, LeBlevec2}.

\begin{figure}
     \centering
\includegraphics[scale=0.6]{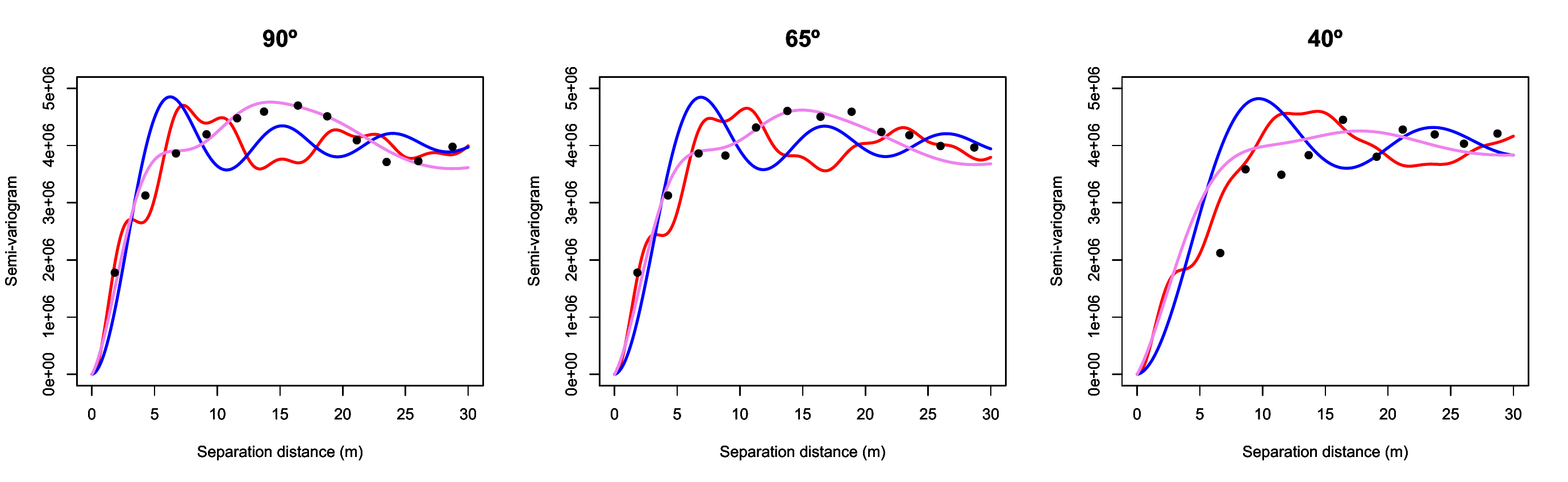}
     \caption{Empirical (black circles) and modeled (solid lines) directional variograms of impedance along directions dipping $90^{\circ}$ (left), $65^{\circ}$ (center) and $40^{\circ}$ (right).  Blue: Model I fitted through CL; Red: Model II fitted through CL; Violet: Model II fitted through LS.}
     \label{fig:variograms}
 \end{figure}

\begin{figure}
    \centering
    \includegraphics[scale=0.27]{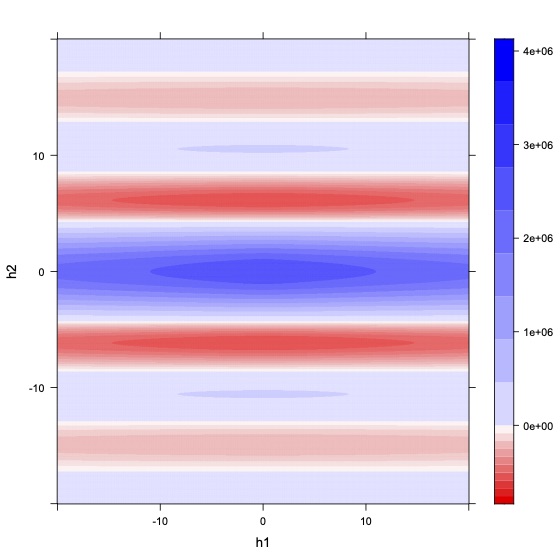}    
     \includegraphics[scale=0.27]{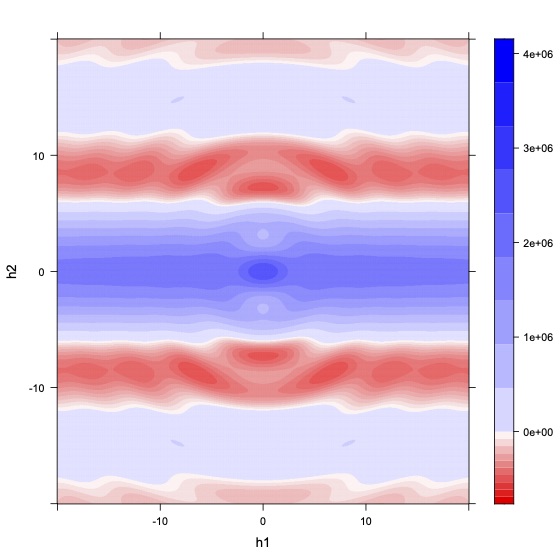}\includegraphics[scale=0.27]{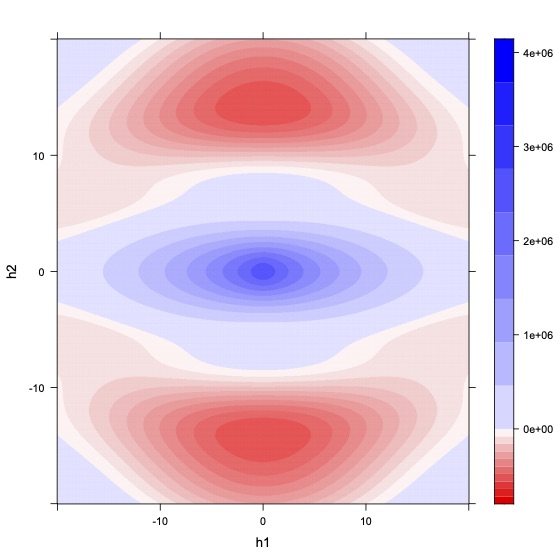}   
    \caption{From left to right: Model I fitted through CL, Model II fitted through CL and Model II fitted through LS.}
    \label{fig:fit_covariances}
\end{figure}

Taking into account this marked axial pattern, characterized by dissimilar scales along the east and depth coordinates, we consider the following models:
\begin{itemize}
    \item {\bf Model I.} A basic construction of the form $$C_{\rm basic}(\bm{h}; \sigma^2,a_1,a_2) = \sigma^2 \exp(-a_1 \|\bm{h}\|) \frac{\sin(a_2 |h_2|)}{a_2 |h_2|},$$
    where $\sigma^2,a_1$ and $a_2$ are positive parameters. 
    \item  {\bf Model II.} We use the previous basic model as a building block and then incorporate derivative information using Proposition \ref{prop3}. The resulting model adopts the form
    $$C(\bm{h}; \sigma^2,a_1,a_2,a_3) = \frac{3\sigma^2}{4}\left[  C_{\rm basic}(\bm{h}; \sigma^2,a_1,a_2) + C_{{\rm derivative}}(\bm{h};a_3) \right],$$
    where 
    $$ C_{{\rm derivative}}(\bm{h};a_3) =  \cos^2(\theta(\bm{h},\bm{u}))\varphi''(\sqrt{a_3}\|\bm{h}\|) +  \sin^2(\theta(\bm{h},\bm{u}))\frac{\varphi'(\sqrt{a_3}\|\bm{h}\|)}{\sqrt{a_3}\|\bm{h}\|},  $$
    with $\bm{u} = [0,1]^\top$ fixed and $\varphi$ of the form \eqref{wave}. Here, $a_3 > 0$ is an additional scale parameter. This model is an example of the variant described in Remark \ref{remark:construction3}.
\end{itemize}
For each model, we estimate the parameters through a composite likelihood (CL) method based on differences \citep{curriero1999composite,10.2307/24309261}. Table \ref{tab:cl} shows the CL estimates together with the value of the objective function at the optimum. For comparison purposes, we also fit a modified version of Model II using an automated least squares (LS) procedure instead of the CL method. For this strategy, we set $\sigma^2 = 4 \times 10^6$, $a_1 = 0.135$, $a_2 = 0.818$ and $a_3=0.067$, in order to obtain a model that matches the structural features of the directional empirical variograms. Figure \ref{fig:variograms} shows the fitted variogram models along three spatial orientations. By construction, Model II that is based on the LS method presents a more accurate description of the empirical variograms, but a poorer log-CL value (Table \ref{tab:cl}). On the contrary,  Models I and II that are based on the CL method do not perfectly match the empirical variograms, a situation that is commonly encountered in practice. To obtain a more comprehensive visualization of the fitted models, Figure \ref{fig:fit_covariances} displays a global plot of the covariances.

\begin{table}
    \centering
    \begin{tabular}{cccccc} \hline\hline 
       Model & $\widehat{\sigma}^2$ & $\widehat{a}_1$ & $\widehat{a}_2$ & $\widehat{a}_3$ & log-CL \\ \hline
       I  & $4.036 \times 10^{6}$  & $1.105 \times 10^{-2}$ &  $0.717$ & $-$ &  $-$13,165,439 \\
       II  & $4.062 \times 10^{6}$  & $3.299 \times 10^{-3}$ & $0.526$ & $2.441$ & $-$13,165,205 \\
     II (based on LS)  & $4 \times 10^{6}$  & $0.135$ & $0.818$ & $0.067$ & $-$13,188,178 \\ 
       \hline \hline
    \end{tabular}
    \caption{Parameters and log-CL of fitted covariance models.}
    \label{tab:cl}
\end{table}

\begin{figure}
    \centering
   \includegraphics[scale=0.4]{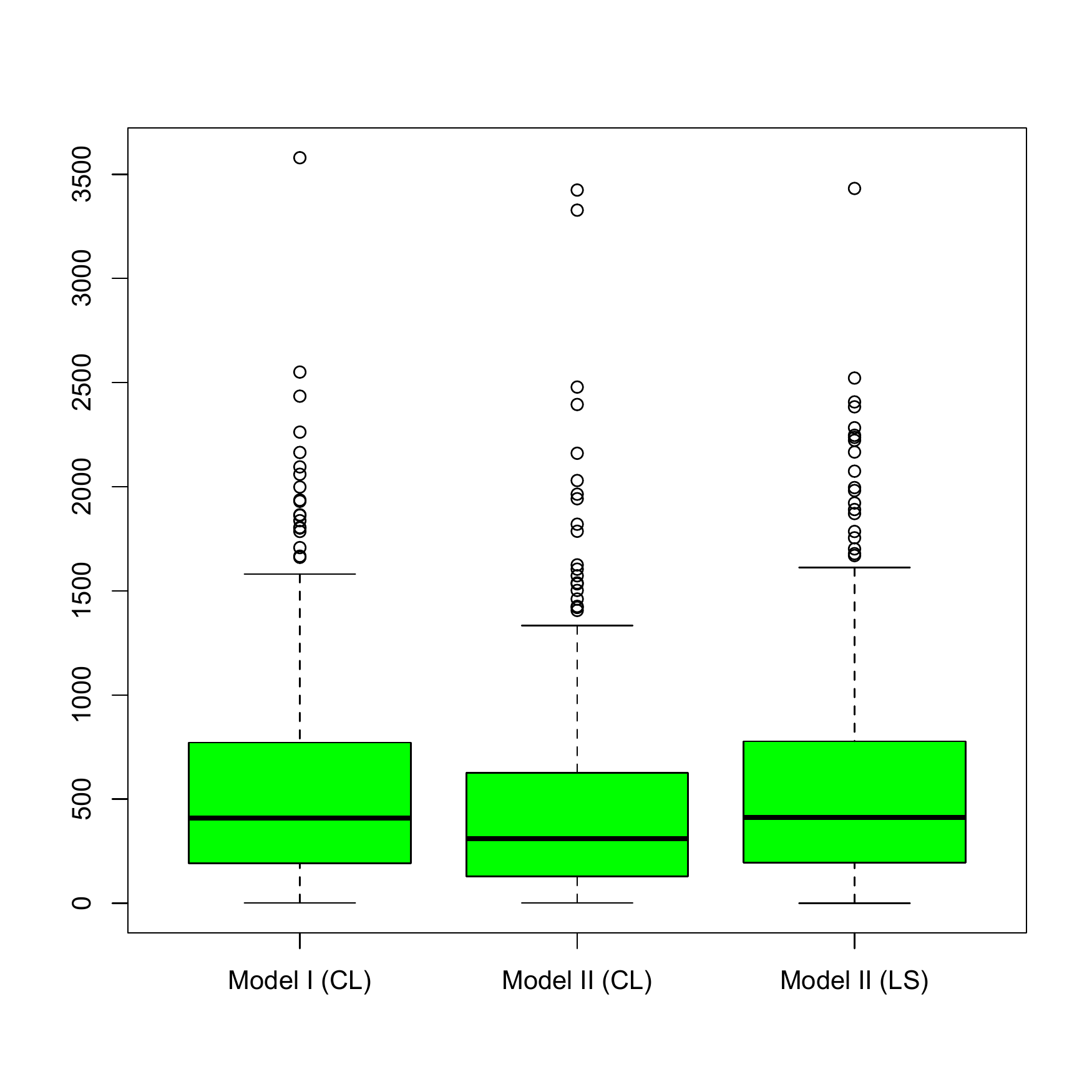}
    \includegraphics[scale=0.4]{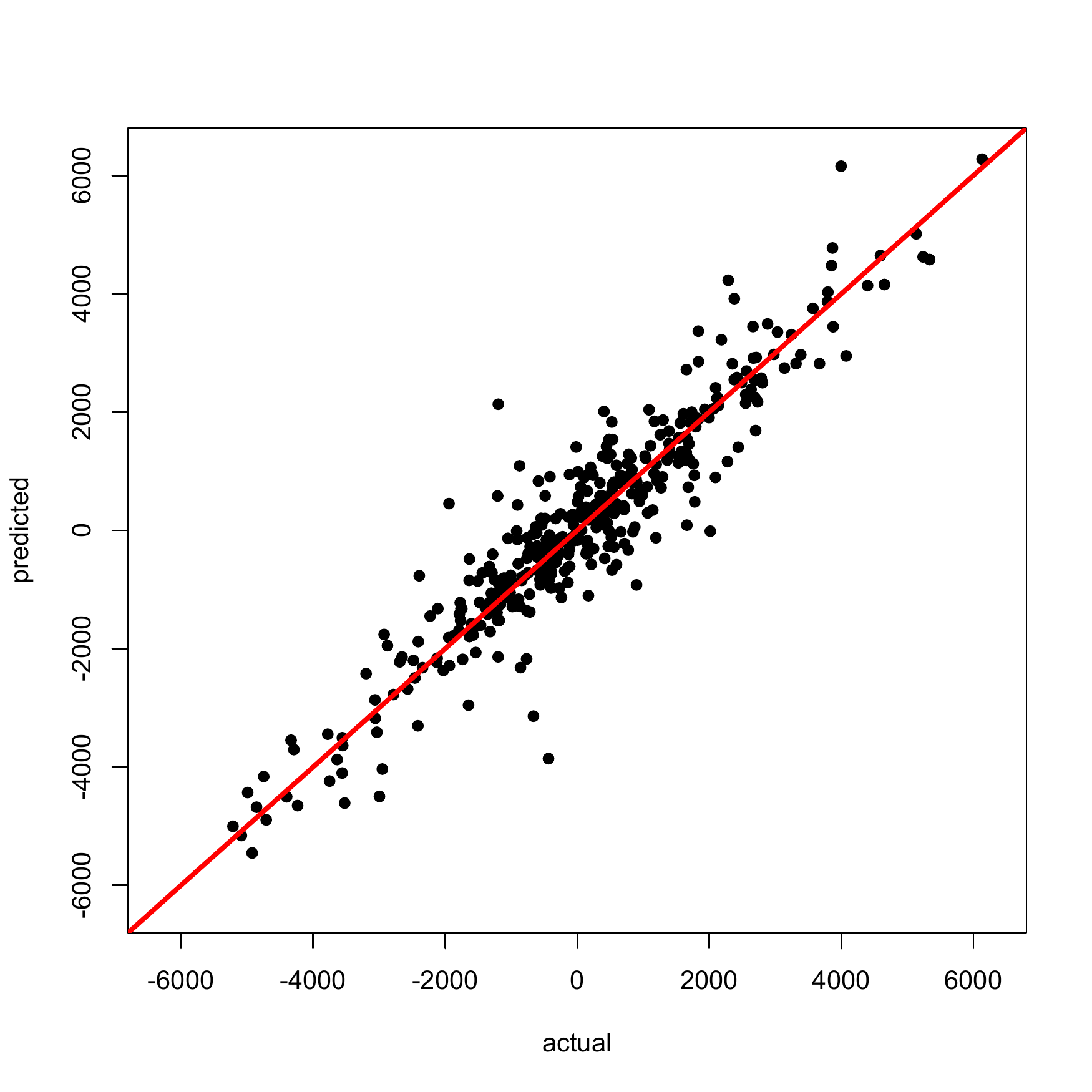}
    \caption{(Left) Comparison of absolute prediction errors among the covariance models. (Right)
 Comparison of actual versus predicted values in the cross-validation study, based on Model II fitted through the CL method.}
    \label{fig:predictions}
\end{figure}

To enhance our analysis, we conduct a split-sample study for model validation, with the $400$ data that have been left out of the model fitting. We apply simple kriging using each model and evaluate the prediction accuracy using metrics such as the root mean square error (RMSE) and mean absolute error (MAE). Among the models that were tested, Model II fitted with the CL method demonstrates a clear advantage, with the RMSE and MAE reduced by $10\%$ to $19\%$ with respect to the other models (see Table \ref{tab:validation}). In Figure \ref{fig:predictions} (left panel), boxplots showing the absolute errors are presented. Model II based on CL outperforms the other models in terms of prediction accuracy. This superiority is evident through noticeably reduced quartiles and upper whisker.  To gain insight into the dispersion of prediction errors, Figure \ref{fig:predictions} (right panel) compares the actual versus predicted values in the validation study, based on Model II fitted through the CL method.

\begin{table}[]
    \centering
    \begin{tabular}{ccc} \hline \hline 
    Model & RMSE & MAE \\ \hline 
      I (based on CL)  &  $741.55$ & $546.12$\\
       II (based on CL)   & $662.27$ & $457.97$\\
      II (based on LS)  &  $767.04$ & $564.76$ \\ \hline \hline 
    \end{tabular}
    \caption{Cross-validation scores: root mean square error (RMSE) and mean absolute error (MAE).}
    \label{tab:validation}
\end{table}

\section{Conclusions}
\label{sec:conclusions}

This work aimed to design new covariance models with complex characteristics. We restricted our attention to models that combine anisotropies and hole effects, and illustrated their practical impact with an application to a geophysical data set. We believe that the pursuit of increasingly flexible models, while maintaining a certain level of simplicity and parsimony, is an area that should continue to be explored. Some recent ideas in this direction can be found in \cite{alegria2021bivariate}, \cite{ma2022beyond}, \cite{10.2307/24311007}  and \cite{BERILD2023100750}, among others.

We illustrated the use of the proposed constructions with well-established families of covariance functions, although our formulations have the potential to be effectively combined with many other parametric families of covariance functions, such as the powered exponential or the hyperbolic models, among others. In particular, employing compactly supported covariances (such as the Gauss hypergeometric covariance) as a starting point provides models that lead to sparse covariance matrices with quite distinctive  structures, allowing for computationally efficient inference \citep{kaufman2008covariance}, prediction \citep{furrer2006covariance} and simulation \citep{Dietrich} techniques. 
 
Extending these results to the multivariate setting, where several coregionalized variables are jointly analyzed and the covariance functions are matrix-valued, presents an interesting area of exploration, albeit accompanied by significant challenges, as the complexity of the models intensifies due to the rapid growth in the number of parameters and the intricate restrictions imposed among them to ensure positive semidefiniteness.

\section*{Acknowledgements}

This work was supported by the National Agency for Research and Development of Chile (ANID), through grants Fondecyt 1210050 (A.A. and X.E.), UTFSM PI$_{-}$LIR$_{-}$23$_{-}$11 (A.A.) and ANID PIA AFB220002 (X.E.).

\bibliography{mybib}
\bibliographystyle{apalike}

\end{document}